\documentclass[11pt]{article}
\usepackage{amsmath,amssymb,geometry,verbatim,algorithm2e,enumerate,float,cite,setspace}
\newtheorem{theorem}{Theorem}
\newtheorem{claim}{Claim}
\newtheorem{lemma}[theorem]{Lemma}
\newtheorem{corollary}[theorem]{Corollary}

\allowdisplaybreaks

\usepackage{setspace}

\begin{document}

\onehalfspace
%\linenumbers

\title{Locally Searching for Large Induced Matchings}

\author{Maximilian F\"{u}rst, Marilena Leichter, Dieter Rautenbach}

\date{}

\maketitle

\begin{center}
Institut f\"{u}r Optimierung und Operations Research, 
Universit\"{a}t Ulm, Ulm, Germany\\
\texttt{maximilian.fuerst@uni-ulm.de}\\
\texttt{marilena.leichter@online.de}\\
\texttt{dieter.rautenbach@uni-ulm.de}\\[3mm]
\end{center}

\begin{abstract}
It is an easy observation that a natural greedy approach yields a 
$\left(d-O(1)\right)$-factor approximation algorithm 
for the maximum induced matching problem in $d$-regular graphs.
The only considerable and non-trivial improvement 
of this approximation ratio was obtained by Gotthilf and Lewenstein 
using a combination of the greedy approach and local search,
where understanding the performance of the local search 
was the challenging part of the analysis.
We study the performance of their local search 
when applied to general graphs,
$C_4$-free graphs,
$\{ C_3,C_4\}$-free graphs,
$C_5$-free graphs, and
claw-free graphs.
As immediate consequences, 
we obtain approximation algorithms 
for the maximum induced matching problem
restricted to the $d$-regular graphs in these classes.

\medskip

\end{abstract}

\section{Introduction}

Finding a maximum induced matching in a given graph is a well-studied difficult problem \cite{stva}, which remains hard even when restricted 
to regular graphs \cite{dumazi}, regular bipartite graphs \cite{dadelo}, or claw-free graphs \cite{koro}.
Next to efficient algorithms for special graph classes \cite{brmo,ca,dumazi,lo},
and fixed parameter tractability \cite{dadelo},
approximation algorithms have been studied.

Duckworth, Manlove, and Zito \cite{dumazi,zi} 
showed that simple and natural greedy strategies yield
$\left(d-O(1)\right)$-factor approximation algorithms
for the maximum induced matching problem 
restricted to $d$-regular graphs.
Up to the constant term, 
this is an immediate consequence 
of the two easy observations fact that 
every induced matching in a $d$-regular graph $G$ contains at most $\frac{m(G)}{2d-1}$ edges, and that 
every maximal induced matching in $G$ contains at least $\frac{m(G)}{2d^2-2d+1}$ edges,
which already implies an approximation ratio of 
$\frac{2d^2-2d+1}{2d-1}=d-\frac{1}{2}+\frac{1}{4d-2}$.
The only considerable and non-trivial improvement 
of this approximation ratio 
was obtained by Gotthilf and Lewenstein \cite{gole} 
who described a 
$\left(0.75d +0.15\right)$-factor approximation algorithm 
for $d$-regular graphs
that combines a greedy approach with local search.
As shown in \cite{ra}, 
for $\{ C_3,C_5\}$-free $d$-regular graphs,
the approximation factor of their method can be improved to $0.708\bar{3}d +0.425$.

Taking a closer look at the algorithm Gotthilf and Lewenstein, 
it turns out that its greedy part {\sc Greedy$(f)$} 
(cf.~Algorithm \ref{alg1} below)
is a kind of preprocessing,
which is rather easy to analyze 
yet ensures an important structural property 
that allows for a better analysis of the local search part {\sc Local Search}
(cf.~Algorithm \ref{alg2} below). 
Their approach can actually improve the approximation ratio by at most a factor of $2$, and is limited by the analysis of {\sc Local Search}.

\bigskip

\begin{algorithm}[H]
{\sc Local Search}\\
\KwIn{A graph $G$.}
\KwOut{An induced matching $M$ of $G$.}
\BlankLine
$M\leftarrow \emptyset$\;
\Repeat{$|M|$ does not increase during one iteration}
{
\If{$M\cup \{ e\}$ is an induced matching of $G$ for some edge $e\in E(G)\setminus M$}
{
$M\leftarrow M\cup \{ e\}$\;
}
\If{$(M\setminus \{ e\})\cup \{ e',e''\}$ is an induced matching of $G$ 
for some three distinct edges $e\in M$ and $e',e''\in E(G)\setminus M$}
{
$M\leftarrow (M\setminus \{ e\})\cup \{ e',e''\}$\;
}
}
\Return $M$\;
\
\caption{The algorithm {\sc Local Search}.}\label{alg2}
\end{algorithm}

\medskip

\noindent Since {\sc Local Search} tries to enlarge the induced matching $M$ by exchanging one edge in $M$ against two other edges in $E(G)\setminus M$, it ensures that the subgraph formed by the so-called {\it private conflict edges} of each individual edge in $M$ is $2K_2$-free;
a structural property that appears quite naturally in the context of induced matchings and the derived strong edge colorings \cite{chgytutr}.

\medskip

\noindent In the present paper we analyze the performance of {\sc Local Search} when applied to general graphs,
$C_4$-free graphs,
$\{ C_3,C_4\}$-free graphs,
$C_5$-free graphs, and
claw-free graphs.
As immediate consequences, 
we obtain approximation algorithms 
for the maximum induced matching problem
restricted to the $d$-regular graphs in these classes.
A byproduct of our results is a much shorter proof of the original result of Gotthilf and Lewenstein \cite{gole}.

Before we proceed to our results, 
we collect some notation and terminology.
We consider only simple, finite, and undirected graphs.
The vertex set and the edge set of a graph $G$
are denoted by $V(G)$ and $E(G)$, respectively.
For a vertex $u$ in a graph $G$,
the neighborhood of $u$ in $G$ is $N_G(u)=\{ v\in V(G):uv\in E(G)\}$,
the closed neighborhood of $u$ in $G$ is $N_G[u]=\{ u\}\cup N_G(u)$, 
and the degree of $u$ in $G$ is $d_G(u)=|N_G(u)|$.
A graph is $d$-regular if all vertices have degree $d$.
For two disjoint sets $X$ and $Y$ of vertices of a graph $G$,
let $E_G(X,Y)$ be the set of edges $uv$ of $G$ with $u\in X$ and $v\in Y$,
and let $m_G(X,Y)=|E_G(X,Y)|$.
For a set $X$ of vertices of $G$, let $E_G(X)$ be the edge set of the subgraph of $G$ induced by $X$, and let $m_G(X)=|E_G(X)|$.

For a set $M$ of edges of a graph $G$, 
let $V(M)$ denote the set of vertices that are incident with an edge in $M$. 
The set $M$ is an induced matching if the subgraph $G[V(M)]$ of $G$ induced by $V(M)$ is $1$-regular.
For a set ${\cal F}$ of graphs, a graph $G$ is ${\cal F}$-free
if it contains no graph in ${\cal F}$ as an induced subgraph.
If ${\cal F}$ contains only one graph $F$, then we write 
$F$-free instead of ${\cal F}$-free.
The cycle of order $n$ is denoted by $C_n$.
The square $G^2$ of a graph $G$ has the same vertex set as $G$, 
and two vertices are adjacent in $G^2$ 
if their distance in $G$ is one or two.
The line graph $L(G)$ of $G$ is the graph whose vertices are the edges of $G$, 
and in which two vertices are adjacent exactly if they share an incident vertex as edges of $G$. 
Note that induced matchings in $G$ correspond to independent sets in $L(G)^2$.

For an edge $e$ of a graph $G$, let 
$$C_G(e)=\{ e\}\cup N_{L(G)^2}(e)=\{ f\in E(G):{\rm dist}_{L(G)}(e,f)\leq 2\},$$ 
and let $c_G(e)=|C_G(e)|$.

For a set $M$ of edges of $G$ and an edge $e$ in $M$, let
$$PC_G(M,e)=C_G(e)\setminus \bigcup_{f\in M\setminus \{ e\}}C_G(f),$$
and let $pc_G(M,e)=|PC_G(M,e)|$.
The set $C_G(e)$ contains $e$ as well as all edges of $G$ that are in {\it conflict} with $e$, that is, that cannot be in an induced matching together with $e$.
The set $PC_G(M,e)$ contains the {\it private conflict} edges of $e$ with respect to some set $M$, which will typically be an induced matching of $G$.

\section{Results}

Throughout this section, let $G$ be a graph of maximum degree $d$ for some $d$ at least $3$,
and let $M$ be an induced matching produced by applying {\sc Local Search} to $G$.

The two essential properties of $M$ are that
\begin{eqnarray}
\mbox{\it for every edge $f$ of $G$, there is some edge $e$ in $M$ with $f\in C_G(e)$}\label{l1}
\end{eqnarray}
and that
\begin{eqnarray}
\mbox{\it $e''\in C_G(e')$ for every edge $e$ in $M$, and every two edges $e'$ and $e''$ in $PC_G(M,e)$.}\label{l2}
\end{eqnarray}
Property (\ref{l1}) means that $M$ is a maximal induced matching, 
and a violation of property (\ref{l2}) would allow {\sc Local Search} to replace $M$ with the larger induced matching $(M\setminus \{ e\})\cup \{ e',e''\}$.

One ingredient of the analysis of {\sc Local Search} is the following upper bound on the number $c_G(xy)$ of edges conflicting with any edge $xy$ of $G$. 
If
$n_{xy}=|N_G(x)\cap N_G(y)|$
and 
$m_{xy}=m_G((N_G(x)\cup N_G(y))\setminus \{ x,y\})$, then
\begin{eqnarray}
c_G(xy) & \leq & 
|\{ xy\}|+d\Big|(N_G(x)\cup N_G(y))\setminus \{ x,y\}\Big|-m_{xy}\nonumber\\
&=& |\{ xy\}|+d\Big(\underbrace{|N_G(x)\cap N_G(y)|}_{=n_{xy}}+
\underbrace{|N_G(x)\setminus N_G[y]|}_{\leq d-1-n_{xy}}+
\underbrace{|N_G(y)\setminus N_G[x]|}_{\leq d-1-n_{xy}}\Big)-m_{xy}\nonumber\\
& \leq & 1+dn_{xy}+2d(d-1-n_{xy})-m_{xy}\nonumber\\
& = & 2d^2-2d+1-\Big(dn_{xy}+m_{xy}\Big)\label{e1}.
\end{eqnarray}
Another ingredient of the analysis of {\sc Local Search} 
are upper bounds on the number of private conflict edges 
for the edges in $M$.
Such bounds are the main concern of the following results.

\medskip

\noindent Let $xy$ be an edge in $M$.

Our first lemma summarizes structural properties of the graph formed by the edges in $PC_G(M,xy)$.

Let 
\begin{itemize}
\item $N_1$ be the set of vertices $u$ in $(N_G(x)\cup N_G(y))\setminus \{ x,y\}$ that are incident with an edge in $PC_G(M,xy)$,
and let 
\item $N_2$ be the set of vertices $u$ in $V(G)\setminus (N_G[x]\cup N_G[y])$ that are incident with an edge in $PC_G(M,xy)$.
\end{itemize}

\begin{lemma}\label{lemma1}
If $G$, $M$, $xy$, $N_1$, and $N_2$ are as above,
then the following statements hold.
\begin{enumerate}[(i)]
\item $PC_G(M,xy)=E_G(\{ x,y\}\cup N_1)\cup E_G(N_1,N_2)$.
\item $N_2$ is independent.
\item If $u\in N_G(x)\cap N_1$ is such that $|N_G(u)\cap N_2|$ is maximum,
then 
$$m_G(N_G(x)\cap N_1,N_2)\leq d-1+(d-3)n_{xu}+m_{ux}.$$
\end{enumerate}
\end{lemma}
{\it Proof:} (i) By the definition of $N_1$ and $N_2$, we have
$PC_G(M,xy)\subseteq E_G(\{ x,y\}\cup N_1)\cup E_G(N_1,N_2)$,
and it remains to show the inverse inclusion.
Since $M$ is an induced matching, we have $xy\in PC_G(M,xy)$.
If $uv$ is an edge in $C_G(xy)$ 
such that $u$ as well as $v$ are both incident with an edge in $PC_G(M,xy)$, then $uv$ also belongs to $PC_G(M,xy)$.
This implies 
$E_G(\{ x,y\}\cup N_1)\cup E_G(N_1,N_2)\subseteq PC_G(M,xy)$,
which completes the proof of (i).

\medskip

\noindent (ii) If $vw$ is an edge between two vertices $v$ and $w$ in $N_2$,
then $vw$ belongs neither to $M$ nor to $C_G(xy)$. 
By (\ref{l1}), 
we obtain $vw\in C_G(e)$ for some edge $e$ in $M$ that is distinct from $xy$,
which implies the contradiction that $v$ or $w$ 
cannot be incident with an edge in $PC_G(M,xy)$. 

\medskip

\noindent (iii) Clearly, $m_G(\{ u\},N_2)\leq d-1-n_{ux}$.
Let $A=(N_G(x)\cap N_1)\cap N_G(u)$ and $B=(N_G(x)\cap N_1)\setminus N_G(u)$. 
Every vertex in $A$ has at most $d-2$ neighbors in $N_2$,
and $|A|\leq n_{ux}$, because every vertex in $A$ is a common neighbor of $x$ and $u$.
By (\ref{l2}),
the choice of $u$ implies that $N_G(v)\cap N_2\subseteq N_G(u)\cap N_2$
for every vertex $v$ in $B$, which implies $m_G(B,N_2)\leq m_{ux}$.
Altogether, we obtain
\begin{eqnarray*}
m_G(N_G(x)\cap N_1,N_2)
&=&m_G(\{ u\},N_2)+m_G(A,N_2)+m_G(B,N_2)\\
&\leq &(d-1-n_{xu})+(d-2)|A|+m_G(B,N_2)\\
&\leq &d-1+(d-3)n_{xu}+m_{ux}. 
\end{eqnarray*}
$\Box$

\medskip

\noindent Our next result establishes a first upper bound on the number $pc_G(M,xy)$ of private conflict edges of $xy$, 
and is slightly more general than the analysis in \cite{gole}.

\begin{theorem}\label{theorem1}
Let $G$ and $M$ be as above.

If $c_G(xy)\geq f$, and $n_{xy}\geq g$ for every edge $xy$ of $G$, then the following statements hold.
\begin{enumerate}[(i)]
\item $\sum\limits_{xy\in M}pc_G(M,xy)\leq 
(4d^2-1-2f-6g)|M|+\sum\limits_{xy\in M}m_{xy}$.
\item $|M|\geq \frac{m(G)}{3d^2-d-f-\left(\frac{d+6}{2}\right)g}$.
\end{enumerate}
\end{theorem}
{\it Proof:} (i) Let $xy$ be an edge in $M$.
Let $N_1$ and $N_2$ be as above.
Let $u$ be as in Lemma \ref{lemma1}(iii).
By (\ref{e1}), we have
$f\leq c_G(ux)\leq 2d^2-2d+1-(dn_{xu}+m_{xu})$, 
which implies $d-1+(d-3)n_{xu}+m_{xu}\leq 2d^2-d-f-3g$.
Now, by symmetry between $x$ and $y$, Lemma \ref{lemma1} implies
\begin{eqnarray*}
pc_G(M,xy) & = & |PC_G(M,xy)|\\
& \leq & |\{ xy\}|+m_G(\{ x,y\},N_1)+m_G(N_1)+m_G(N_G(x)\cap N_1,N_2)+m_G(N_G(y)\cap N_1,N_2)\\
& \leq & 1+(2d-2)+m_{xy}+(2d^2-d-f-3g)+(2d^2-d-f-3g)\\
& = & 4d^2-1-2f-6g+m_{xy}.
\end{eqnarray*}
Adding this inequality over all edges in $M$ yields the desired bound.

\medskip

\noindent (ii) If $p$ is the number of pairs $(e,f)$ with $e\in M$ and $f\in C_G(e)$, then
\begin{eqnarray*}
p & \leq & \sum\limits_{xy\in M}c_G(xy)
\stackrel{(\ref{e1})}{\leq} (2d^2-2d+1-dg)|M|-\sum\limits_{xy\in M}m_{xy}.
\end{eqnarray*}
Furthermore, since the only edges $f$ of $G$, 
for which there is only one edge $e$ in $M$ with $f\in C_G(e)$, 
are those in $\bigcup\limits_{xy\in M}PC_G(M,xy)$, we obtain
\begin{eqnarray*}
p & \geq & 2m(G)-\sum\limits_{xy\in M}pc_G(M,xy)
\stackrel{(i)}{\geq} 2m(G)-(4d^2-1-2f-6g)|M|-\sum\limits_{xy\in M}m_{xy}.
\end{eqnarray*}
Combining the upper and lower bounds on $p$ yields the desired bound.
$\Box$

\medskip

\noindent With Theorem \ref{theorem1} at hand, 
it is now easy to recover the result of Gotthilf and Lewenstein \cite{gole}.

\begin{corollary}[Gotthilf and Lewenstein \cite{gole}]\label{corollary1}
There is a polynomial time 
$\left(\frac{3}{4}d+\frac{1}{8}+\frac{1}{16d-8}\right)$-factor approximation algorithm 
for the maximum induced matching problem 
in $d$-regular graphs.
\end{corollary}
{\it Proof:} Let $G$ be a given $d$-regular graph.
Applying {\sc Greedy$(f)$} to $G$ with $f=\frac{3d^2-d}{2}$ yields an output $(M,G')$ with 
$|M|\geq \frac{m(G)-m(G')}{f}$, and 
$c_{G'}(xy)\geq f$ for every edge $xy$ of $G'$.
Now, by Theorem \ref{theorem1}(ii), 
applying {\sc Local Search} to $G'$ yields an induced matching 
$M'$ of $G'$ with 
$|M'|\geq \frac{m(G')}{3d^2-d-f}=\frac{m(G')}{f}$,
where the last equality follows from the choice of $f$.
It is easy to see (cf.~\cite{ra}) that $M\cup M'$ is an induced matching of $G$
with $|M\cup M'|\geq \frac{m(G)}{f}$.
Since $\nu_s(G)\leq \frac{m(G)}{2d-1}$, 
we obtain $\frac{\nu_s(G)}{|M\cup M'|}\leq \frac{f}{2d-1}=\frac{3}{4}d+\frac{1}{8}+\frac{1}{16d-8}$,
which completes the proof. $\Box$

\medskip

\begin{algorithm}[H]
{\sc Greedy$(f)$}\\
\KwIn{A graph $G$.}
\KwOut{A pair $(M,G')$ such that $M$ is an induced matching of $G$, and $G'$ is a subgraph of $G$.}
\BlankLine
$M\leftarrow \emptyset$;
$G_0\leftarrow G$; 
$i\leftarrow 1$\;
\While{$\min\{ c_{G_{i-1}}(e):e\in E(G_{i-1})\}\leq f$}
{
Choose an edge $e_i$ of $G_{i-1}$ with $c_{G_{i-1}}(e_i)\leq f$\;
$M\leftarrow M\cup \{ e_i\}$;
$G_i\leftarrow G_{i-1}-C_{G_{i-1}}(e_i)$;
$i\leftarrow i+1$;
}
\Return $(M,G_{i-1})$\;
\
\caption{The algorithm {\sc Greedy$(f)$}.}\label{alg1}
\end{algorithm}

\medskip

\noindent Note that first applying {\sc Greedy$(f)$} corresponds to a preprocessing ensuring a lower bound on $\min\{ c_G(xy):xy\in E(G)\}$,
which is important for the performance of {\sc Local Search}.

Our next two results concern the performance of {\sc Local Search}
for $C_4$-free and $\{ C_3,C_4\}$-free graphs.
Both rely on the two step approach from Theorem \ref{theorem1}
of upper bounding $pc_G(M,xy)$, and double-counting $p$.

\begin{theorem}\label{theorem2}
Let $G$ and $M$ be as above.

If $G$ is $C_4$-free, then the following statements hold.
\begin{enumerate}[(i)]
\item If $xy$ is an edge in $M$, and $N_1$ and $N_2$ are as above, then 
$m_G(N_1,N_2)\leq \max\left\{ d,\frac{d^2}{4}\right\}$.
\item $|M|\geq \frac{m(G)}{\max\left\{d^2+\frac{d}{2},\frac{9}{8}d^2\right\}}$.
\end{enumerate}
\end{theorem}
{\it Proof:} (i) 
Let $N_1'$ be the set of those $u$ in $N_1$ that have a neighbor in $N_2$.
Let 
$N_x=N'_1\setminus N_G(y)$,
$N_y=N'_1\setminus N_G(x)$, and
$N_{xy}=N'_1\cap N_G(x)\cap N_G(y)$.

Suppose, for a contradiction, 
that $u$ and $u'$ are two non-adjacent vertices in $N_x\cup N_{xy}$.
Let $v$ be a neighbor of $u$ in $N_2$,
and let $v'$ be a neighbor of $u'$ in $N_2$.
Since $G$ is $C_4$-free, the vertices $u$ and $u'$ have no common neighbor in $N_2$; 
in particular, the vertices $v$ and $v'$ are distinct.
Since, by Lemma \ref{lemma1}(ii), the vertices $v$ and $v'$ are non-adjacent, we obtain a contradiction to (\ref{l2}).
Hence, $N_x\cup N_{xy}$, and, by symmetry, also $N_y\cup N_{xy}$ are cliques in $G$.

Suppose, for a contradiction,
that the two vertices $u$ in $N_x$ and $u'$ in $N_y$ both have at least two neighbors in $N_2$.
Since $G$ is $C_4$-free, the vertices $u$ and $u'$ are non-adjacent.
By (\ref{l2}), and since, by Lemma \ref{lemma1}(ii), the set $N_2$ is independent,
it follows that $u$ and $u'$ have two common neighbors in $N_2$,
contradicting the $C_4$-freeness of $G$.
Hence, we may assume, by symmetry between $x$ and $y$, 
that each vertex in $N_y$ has exactly one neighbor in $N_2$.

It follows that, 
if $d_x=|N_x|$, $d_{xy}=|N_{xy}|$, and $d_y=|N_{xy}|$, then 
\begin{itemize}
\item every vertex in $N_x$ has at most $(d-1)-(d_x+d_{xy}-1)$ neighbors in $N_2$,
\item every vertex in $N_{xy}$ has at most $(d-2)-(d_x+d_{xy}+d_y-1)$ neighbors in $N_2$, and
\item every vertex in $N_y$ has one neighbor in $N_2$.
\end{itemize}
Note that 
$d_x+d_{xy}\leq d_G(x)-1\leq d-1$ and
$d_y+d_{xy}\leq d_G(y)-1\leq d-1$, and that 
\begin{eqnarray}
\mbox{\it the function $x\mapsto (d-x)x$ has its unique maximum $d^2/4$ at $x=d/2$}.\label{c1}
\end{eqnarray}
First, we assume that $d_x$ and $d_y$ are both positive.
If $u\in N_x$ and $v\in N_y$, 
then, by (\ref{l2}), $yv\in C_G(e)$ for every edge $e$ in $PC_G(M,xy)$ incident with $u$,
which implies that $v$ is adjacent to every neighbor of $u$ in $N_2$.
Since $v$ has exactly one neighbor, say $w$, in $N_2$, 
and the choice of $u$ within $N_x$ was unrestricted, 
it follows that $N_G(u)\cap N_2=\{ w\}$ for every vertex $u$ in $N_x$.
Now, by (recovered) symmetry between $x$ and $y$, it follows that $N_G(v)\cap N_2=\{ w\}$ for every vertex $v$ in $N_y$.

If $d_x+d_{xy}+d_y\geq d-1$, then 
\begin{eqnarray*}
m_G(N_1,N_2) & \leq & d_x+\underbrace{\Big((d-2)-(d_x+d_{xy}+d_y-1)\Big)}_{\leq 0}d_{xy}+d_y
\leq d_x+d_y \leq d_G(w) \leq d.
\end{eqnarray*}
If $d_x+d_{xy}+d_y<d-1$, then 
\begin{eqnarray*}
m_G(N_1,N_2) & \leq & d_x+\underbrace{\Big((d-2)-(d_x+d_{xy}+d_y-1)\Big)}_{\geq 1}
d_{xy}+d_y\\
&\leq & \Big(d-1-(d_x+d_{xy}+d_y)\Big)(d_x+d_{xy}+d_y)\\
& \stackrel{(\ref{c1})}{\leq} & \frac{(d-1)^2}{4}.
\end{eqnarray*}
Next, we assume that $d_y=0$.
In this case, we obtain
\begin{eqnarray*}
m_G(N_1,N_2) & \leq & 
\Big((d-1)-(d_x+d_{xy}-1)\Big)d_x+\Big((d-2)-(d_x+d_{xy}+d_y-1)\Big)d_{xy}\\
& \leq & 
\Big(d-(d_x+d_{xy})\Big)(d_x+d_{xy})\\
& \stackrel{(\ref{c1})}{\leq} & 
\frac{d^2}{4}.
\end{eqnarray*}
Finally, if $d_x=0$, then a similar argument implies $m_G(N_1,N_2) \leq \frac{d^2}{4}$,
which completes the proof of (i).

\medskip

\noindent (ii) By (i), we obtain, similarly as in the proof of Theorem \ref{theorem1}(i),
$$pc_G(M,xy)\leq 1+(2d-2)+m_{xy}+m_G(N_1,N_2)\stackrel{(i)}{\leq} 
\max\left\{ 3d-1,\frac{d^2}{4}+2d-1\right\}+m_{xy},$$
and, hence,
$$\sum\limits_{xy\in M}pc_G(M,xy)\leq  
\max\left\{ 3d-1,\frac{d^2}{4}+2d-1\right\}|M|+\sum\limits_{xy\in M}m_{xy}$$
Counting the $p$ pairs as in the proof of Theorem \ref{theorem1}(ii), 
we obtain
$p\leq (2d^2-2d+1)|M|-\sum\limits_{xy\in M}m_{xy}$
and 
$p\geq 2m(G)-\max\left\{ 3d-1,\frac{d^2}{4}+2d-1\right\}|M|-\sum\limits_{xy\in M}m_{xy}$,
which yields the desired bound.
$\Box$

\begin{theorem}\label{theorem3}
Let $G$ and $M$ be as above.

If $G$ is $\{ C_3,C_4\}$-free, then the following statements hold.
\begin{enumerate}[(i)]
\item If $xy$ is an edge in $M$, then $pc_G(M,xy)\leq 2d-1$.
\item $|M|\geq \frac{m(G)}{d^2}$.
\end{enumerate}
\end{theorem}
{\it Proof:} (i) Let $N_1$ and $N_2$ be as above.
Since $G$ is $C_3$-free, the set $N_1$ partitions into the two independent sets 
$N_x=N_1\cap N_G(x)$ and $N_y=N_1\cap N_G(y)$.
If one of the two sets $N_x$ or $N_y$ contains two vertices with a neighbor in $N_2$,
then we obtain a similar contradiction as in the proof of Theorem \ref{theorem2}(i).
Hence, both sets contain at most one such vertex.

If both sets $N_x$ and $N_y$ contain no vertex with a neighbor in $N_2$, 
then $pc_G(M,xy)\leq 2d-1$.
If $N_x$ contains a vertex with a neighbor in $N_2$ but $N_y$ does not, then $N_y$ is empty, and, hence, $pc_G(M,xy)\leq 2d-1$.
Finally, if the vertex $u$ in $N_x$ as well as the vertex $v$ in $N_y$ 
both have a neighbor in $N_2$, then, since, by (\ref{l2}), 
$xu\in C_G(e)$ for every edge $e$ in $PC_G(M,xy)$ that is incident with $v$,
and 
$yv\in C_G(f)$ for every edge $f$ in $PC_G(M,xy)$ that is incident with $u$,
it follows that 
the vertices $u$ and $v$ both have exactly the same neighbors in $N_2$.
Since $G$ is $\{ C_3,C_4\}$-free, the vertices $u$ and $v$ have a unique neighbor $w$ in $N_2$,
in particular, $N_2=\{ w\}$.
If $N_x$ contains a vertex $u'$ distinct from $u$, then $xu'\not\in C_G(vw)$,
which is a contradiction to (\ref{l2}).
Hence, $N_x$, and, by symmetry, also $N_y$ both contain exactly one vertex, 
and $pc_G(M,xy)=5\leq 2d-1$.

\medskip

\noindent (ii) If $p$ is as in the previous proofs, we obtain
$p\leq (2d^2-2d+1)|M|$ and $p\geq 2m(G)-(2d-1)|M|$,
which yields the desired bound. $\Box$

\medskip

\noindent Since Theorem \ref{theorem2} and Theorem \ref{theorem3}
do not require a lower bound on $\max\{ c_G(xy):xy\in E(G)\}$,
they imply that {\sc Local Search} alone is an approximation algorithm
in the regular case.

\begin{corollary}\label{corollary2}
{\sc Local Search} is 
\begin{enumerate}[(i)]
\item a polynomial time 
$\left(\frac{9}{16}d+\frac{33}{80}\right)$-factor approximation algorithm 
for the maximum induced matching problem in $d$-regular $C_4$-free graphs,
and
\item a polynomial time 
$\left(\frac{d}{2}+\frac{1}{4}+\frac{1}{8d-4}\right)$-factor approximation algorithm 
for the maximum induced matching problem in $d$-regular $\{ C_3,C_4\}$-free graphs.
\end{enumerate}
\end{corollary}
{\it Proof:} Using Theorem \ref{theorem2}(ii), 
part (i) follows as in Corollary \ref{corollary1}, 
since $\frac{\max\left\{d^2+\frac{d}{2},\frac{9}{8}d^2\right\}}{2d-1}\leq \frac{9}{16}d+\frac{33}{80}$ for $d\geq 3$.
Similarly, part (ii) follows from Theorem \ref{theorem3}(ii).
$\Box$

\medskip

\noindent Next, we consider $C_5$-free graphs.
As it turns out, a better bound can be obtained for these graphs,
if the lower and upper bounds on $p$ are not handled separately.

\begin{theorem}\label{theorem4}
Let $G$ and $M$ be as above.

If $G$ is $C_5$-free, then the following statements hold.
\begin{enumerate}[(i)]
\item If $xy$ is an edge in $M$, and $N_1$ and $N_2$ are as above, then
$$m_G(N_1,N_2)\leq (d-1-d_y)d_x+(d-2)d_{xy}+(d-1-d_x)d_y$$
for non-negative integers $d_x,d_{xy}$, and $d_y$ with 
$d_x+d_{xy}\leq d-1$,
$d_y+d_{xy}\leq d-1$, and
$d_{xy}\leq n_{xy}$.
\item $|M|\geq \frac{m(G)}{\frac{3}{2}d^2-d+\frac{1}{2}}$.
\end{enumerate}
\end{theorem}
{\it Proof:} (i) Let $N_x$, $N_y$, $N_{xy}$, $d_x$, $d_y$, and $d_{xy}$ be as in the proof of Theorem \ref{theorem2}(i).
Clearly, $d_x+d_{xy}\leq d_G(x)-1\leq d-1$ and $d_y+d_{xy}\leq d_G(y)-1\leq d-1$. Furthermore, $n_{xy}=|N_G(x)\cap N_G(y)|\geq d_{xy}$.

Since $G$ is $C_5$-free, (\ref{l2}) implies that 
every vertex in $N_x$ is adjacent to every vertex in $N_y$, that is, 
every vertex in $N_x$ has at most $d-1-d_y$ neighbors in $N_2$,
and 
every vertex in $N_y$ has at most $d-1-d_x$ neighbors in $N_2$.
Since every vertex in $N_{xy}$ has at most $d-2$ neighbors in $N_2$,
part (i) follows.

\medskip

\noindent (ii) For every edge $xy$ in $M$, (i) implies
\begin{eqnarray}
pc_G(M,xy) & \leq & 1+(2d-2)+m_{xy}+\Big((d-1-d_y)d_x+(d-2)d_{xy}+(d-1-d_x)d_y\Big),\label{e2}
\end{eqnarray}
where $d_x$, $d_y$, and $d_{xy}$ are non-negative integers that depend on $xy$, and satisfy the restrictions stated in (i).

If $p$ is as in the previous proofs, then 
$$2m(G)-\sum\limits_{xy\in M}pc_G(M,xy)\leq p\leq \sum\limits_{xy\in M}(2d^2-2d+1-dd_{xy}-m_{xy}),$$
where the second inequality follows from (\ref{e1}) using $d_{xy}\leq n_{xy}$.
By (\ref{e2}), we obtain
\begin{eqnarray}
\sum\limits_{xy\in M}
\Big(2d^2-dd_{xy}+(d-1-d_y)d_x+(d-2)d_{xy}+(d-1-d_x)d_y\Big)\geq 2m(G).\label{e3}
\end{eqnarray}
Since
$\max\{ (d-1-y)x+(d-1-x)y:0\leq x,y\leq d-1\}=(d-1)^2$,
we obtain, 
using (\ref{e3}) and $d_{xy}\geq 0$, 
that $(2d^2+(d-1)^2)|M|\geq 2m(G)$,
which completes the proof of (ii).
$\Box$

\medskip

\noindent For regular graphs, we obtain an approximation guarantee as above.

\begin{corollary}\label{corollary3}
{\sc Local Search} is a polynomial time 
$\left(\frac{3}{4}d-\frac{1}{8}+\frac{3}{16d-8}\right)$-factor approximation algorithm 
for the maximum induced matching problem in $d$-regular $C_5$-free graphs.
\end{corollary}
{\it Proof:} This follows by evaluating $\frac{\frac{3}{2}d^2-d+\frac{1}{2}}{2d-1}$. $\Box$

\medskip

\noindent Now, we consider claw-free graphs.
Our next result illustrates the consequences of applying the above arguments to these graphs.

\begin{lemma}\label{lemma2}
If $G$ and $M$ are as above, and $G$ is $K_{1,3}$-free, then $|M|\geq \frac{m(G)}{d^2+d-1}$.
\end{lemma}
{\it Proof:} If $xy$ is an edge in $M$, and $N_1$ and $N_2$ are as above,
then claw-freeness and Lemma \ref{lemma1}(ii) together imply 
that every vertex in $N_1$ has at most one neighbor in $N_2$,
which implies $m_G(N_1,N_2)\leq |N_1|\leq 2d-2$, 
and, hence,
$pc_G(M,xy)\leq 1+(2d-2)+m_{xy}+(2d-2)=4d-3+m_{xy}$.
Double-counting $p$ as before yields the desired result. $\Box$

\medskip

\noindent In contrast to the classes of graphs that we considered before, 
no edge $xy$ in a claw-free graph $G$ of maximum degree $d$ 
can have up to $2d^2+O(d)$ many conflict edges,
that is, claw-freeness implies a better upper bound on $c_G(xy)$,
which implies that already maximal induced matchings cannot be too small.
The next result quantifies this observation.

Recall that two vertices $u$ and $v$ are {\it true twins} in a graph $G$
if $N_G[u]=N_G[v]$.

\begin{theorem}\label{theorem5}
If $G$ is a claw-free graph of maximum degree at most $d$ for some $d$ at least $3$,
and $M$ is a maximal induced matching, then
$|M|\geq \frac{m(G)}{\frac{7}{6}d^2+d}.$
\end{theorem}
{\it Proof:} Let $xy$ be an edge in $M$.
Let 
$N_x=N_G(x)\setminus N_G[y]$,
$N_y=N_G(y)\setminus N_G[x]$, and
$N_{xy}=N_G(x)\cap N_G(y)$, 
and let 
$d_x=|N_x|$,
$d_y=|N_y|$, and
$d_{xy}=|N_{xy}|$.

Since $c_G(xy)\leq 1+d(d_x+d_{xy}+d_y)-m_{xy}$, 
in order to upper bound $c_G(xy)$, 
we need a lower bound on $m_{xy}$.
Let $H=G[N_x\cup N_{xy}\cup N_y]$.
By definition, $m_{xy}=m(H)$.

Since $G$ is claw-free,
\begin{enumerate}[(i)]
\item $N_x$ and $N_y$ are cliques in $H$, and
\item $\alpha(H[N_x\cup N_{xy}])\leq 2$ and $\alpha(H[N_{xy}\cup N_y])\leq 2$.
\end{enumerate}
\begin{claim}\label{claim1}
$m(H) \geq 
{d_x\choose 2}
+
{p\choose 2}
+
{d_{xy}-p\choose 2}
+
{d_y\choose 2}
+p(d_x+d_y)$
for some non-negative integer $p\leq \frac{d_{xy}}{2}$.
\end{claim}
{\it Proof of Claim \ref{claim1}:}
Suppose that, subject to (i) and (ii), 
the graph $H$ is such that $m(H)$ is minimum. 
Our first goal is to show that $H[N_{xy}]$ is the disjoint union of cliques.
Suppose, for a contradiction, 
that $uvw$ is an induced path of order three in $H[N_{xy}]$.
If $d_H(u)<d_H(v)$, then removing $v$, and introducing a new vertex that is a true twin of $u$ yields a graph with less edges than $H$
that satisfies (i) and (ii),
which is a contradiction.
If $d_H(u),d_H(w)\geq d_H(v)$, then removing $u$ and $w$, and introducing two new vertices
that are true twins of $v$ yields a graph with less edges than $H$
that satisfies (i) and (ii),
which is a contradiction.
Hence, by (ii), the graph $H[N_{xy}]$ is the union of two cliques,
$C_1$ of order $p$ and $C_2$ of order $d_{xy}-p$,
where $0\leq p\leq \frac{d_{xy}}{2}$.
Note that $p=0$ corresponds to the case that $G[N_{xy}]$ is just one clique.
By (ii), every vertex in $N_x$ is adjacent to all vertices in $C_1$ or to all vertices in $C_2$. By the choice of $H$, and symmetry between $N_x$ and $N_y$, we obtain, that $N_x\cup N_y$ is completely joined to $C_1$, and the desired lower bound on $m(H)$ follows.
$\Box$

\medskip

\noindent Let
$\bar{d}_x=\frac{d_x}{d}$,
$\bar{d}_{xy}=\frac{d_{xy}}{d}$,
$\bar{d}_y=\frac{d_y}{d}$, and
$\bar{p}=\frac{p}{d}$.
Since $\frac{1}{d^2}+\frac{d_x+p+(d_{xy}-p)+d_y}{2d^2}\leq \frac{1}{d^2}+\frac{2d-2}{2d^2}=\frac{1}{d}$,
Claim \ref{claim1} implies
\begin{eqnarray}
c_G(xy) 
& \leq & 1+d(d_x+d_{xy}+d_y)-m(H)\nonumber\\
& \leq & 1+d(d_x+d_{xy}+d_y)-{d_x\choose 2}
-
{p\choose 2}
-
{d_{xy}-p\choose 2}
-
{d_y\choose 2}
-p(d_x+d_y)\nonumber\\
& \leq & d^2\Bigg(
\underbrace{\bar{d}_x+\bar{d}_{xy}+\bar{d}_y
-\frac{1}{2}\bar{d}_x^2
-\frac{1}{2}\bar{p}^2
-\frac{1}{2}\left(\bar{d}_{xy}-\bar{p}\right)^2
-\frac{1}{2}\bar{d}_y^2
-\bar{p}(\bar{d}_x+\bar{d}_y)}_{=:f\left(\bar{d}_x,\bar{d}_{xy},\bar{d}_y,\bar{p}\right)}
+\frac{1}{d}
\Bigg),\label{e4}
\end{eqnarray}
where 
\begin{eqnarray}
\bar{d}_x+\bar{d}_{xy} & \leq & 1,\label{ec1}\\
\bar{d}_y+\bar{d}_{xy} & \leq & 1,\label{ec2}\\
\bar{p} & \leq & \frac{\bar{d}_{xy}}{2},\mbox{ and}\label{ec3}\\
\bar{d}_x,\bar{d}_{xy},\bar{d}_y,\bar{p} & \geq & 0.\label{ec4}
\end{eqnarray}
Suppose that,
subject to (\ref{ec1}) to (\ref{ec4}),
the values of $\bar{d}_x$, $\bar{d}_{xy}$, $\bar{d}_y$, and $\bar{p}$
are chosen such that $f\left(\bar{d}_x,\bar{d}_{xy},\bar{d}_y,\bar{p}\right)$ assumes its maximum value $f_{\max}$.
Since $f\left(\bar{d}_x,\bar{d}_{xy},\bar{d}_y,\bar{p}\right)$
equals $-\frac{1}{2}\bar{d}_x^2+(1-\bar{p})\bar{d}_x$
plus a function that does not depend on $\bar{d}_x$, 
and, hence,
$\frac{\partial}{\partial \bar{d}_x}f\left(\bar{d}_x,\bar{d}_{xy},\bar{d}_y,\bar{p}\right)=1-\bar{p}-\bar{d}_x$,
it follows that condition (\ref{ec1}), and, by symmetry, 
also condition (\ref{ec2}) holds with equality.
Setting $g\left(\bar{d}_{xy},\bar{p}\right):=f\left(1-\bar{d}_{xy},\bar{d}_{xy},1-\bar{d}_{xy},\bar{p}\right)$, we obtain
$$g\left(\bar{d}_{xy},\bar{p}\right)=
1-\frac{3}{2}\bar{d}_{xy}^2-\bar{p}^2+\bar{d}_{xy}-2\bar{p}+3\bar{d}_{xy}\bar{p}.$$
First, suppose that $\bar{d}_{xy}\leq \frac{2}{3}$.
Since 
$\frac{\partial}{\partial \bar{d}_{xy}}g\left(\bar{d}_{xy},\bar{p}\right)=1+3\bar{p}-3\bar{d}_{xy}$,
we obtain that $\bar{d}_{xy}=\bar{p}+\frac{1}{3}$,
and that $0\leq \bar{p}\leq \frac{1}{3}$.
Note that in this case
\begin{eqnarray}
g\left(\bar{p}+\frac{1}{3},\bar{p}\right)=\frac{7}{6}-\bar{p}+\frac{1}{2}\bar{p}^2\leq \frac{7}{6}.\label{fmax1}
\end{eqnarray}
Next, suppose that $\bar{d}_{xy}\geq \frac{2}{3}$.
Since 
$\frac{\partial}{\partial \bar{p}}g\left(\bar{d}_{xy},\bar{p}\right)=
3\bar{d}_{xy}-2-2\bar{p}$,
we obtain that $\bar{p}=\frac{3}{2}\bar{d}_{xy}-1$,
and that $\frac{2}{3}\leq\bar{d}_{xy}\leq 1$.
Note that in this case
\begin{eqnarray}
g\left(\bar{d}_{xy},\frac{3}{2}\bar{d}_{xy}-1\right)=
2-2\bar{d}_{xy}+\frac{3}{4}\bar{d}_{xy}^2\leq 1.\label{fmax2}
\end{eqnarray}
In view of (\ref{fmax1}) and (\ref{fmax2}), we obtain $f_{\max}=\frac{7}{6}$, which, by (\ref{e4}), implies
$c_G(xy) \leq \frac{7}{6}d^2+d$.
If $p$ is as in the previous proofs, then
$p\leq \sum\limits_{xy\in M}c_G(xy)\leq \left(\frac{7}{6}d^2+d\right)|M|$,
and, by the maximality of $M$, we obtain $p\geq m(G)$,
which implies the desired lower bound on $|M|$. $\Box$

\medskip

\noindent We close with the obvious corollary.

\begin{corollary}\label{corollary4}
{\sc Local Search} is a 
polynomial time 
$\left(\frac{1}{2}d+\frac{3}{4}-\frac{1}{8d-4}\right)$-factor approximation algorithm 
for the maximum induced matching problem in claw-free $d$-regular graphs.
Furthermore, choosing any maximal induced matching is a 
polynomial time 
$\left(\frac{7}{12}d+\frac{19}{24}+\frac{19}{48d-24}\right)$-factor approximation algorithm 
for the same problem. 
\end{corollary}

\medskip

\noindent It seems possible to generalize Theorem \ref{theorem5} and Corollary \ref{corollary4} to $K_{1,r}$-free graphs. 
Combining {\sc Local Search} with some greedy preprocessing might further improve the approximation guarantees obtained in this paper.

\pagebreak

\end{document}